\documentclass[final,times,twocolumn]{elsarticle}

\usepackage{lineno}

\usepackage{amssymb}

\usepackage{blindtext} 
\usepackage{booktabs}
\usepackage{eurosym}
\usepackage{siunitx}
\usepackage{chemformula}
\DeclareSIUnit{\EUR}{\mbox{\euro}}
\DeclareSIUnit{\USD}{\mbox{\$}}
\DeclareSIUnit{\year}{\mbox{y}}
\DeclareSIUnit{\tCO}{\mbox{t\ch{CO2}}}
\DeclareSIUnit{\betaCo}{\mbox{m$^{2\,\beta}$}}
\DeclareSIUnit{\wtpercent}{wt\%}
\DeclareSIUnit{\ct}{ct}
\usepackage{mathrsfs}
\usepackage{float}
\usepackage{threeparttable}
\usepackage{tabularx}
\usepackage{amsmath}

\usepackage{nomencl}
\usepackage{etoolbox}
\usepackage{svg}
\usepackage{graphicx}
\usepackage{float}
\usepackage{comment}
\usepackage{multirow}
\usepackage[justification=centering]{caption}
\usepackage[super]{nth}
\usepackage{xurl}
\usepackage{soul}

\usepackage{verbatim}

\usepackage{xargs}    

\usepackage{regexpatch}
\makeatletter
\regexpatchcmd\ps@pprintTitle
  {\cE\}\ \cE\}}{\cE\}\cE\}}
  {}{\FailedToPatch}
\makeatother

\begin{document}

\begin{frontmatter}



\title{Evaluating synthetic fuel production: A case study on the influence of electricity and \ch{CO2} price variations}

\journal{Case Studies in Thermal Engineering}

\author[inst1]{David Huber\corref{corrauthor}}\ead{david.huber@tuwien.ac.at}
\author[inst1]{Felix Birkelbach}
\author[inst1]{René Hofmann}

\cortext[corrauthor]{Corresponding author}

\address[inst1]{TU Wien, Institute of Energy Systems and Thermodynamics, Getreidemarkt 9/BA, 1060 Vienna, Austria}

\begin{abstract}
To combat climate change, we need to reduce emissions from the transport sector. Synthetic fuels are a long-term solution for aviation, maritime and heavy machinery. Large-scale use requires cost-effectiveness, efficient production and resilience to price changes. In this case study, we simultaneously optimize the cell voltage of the solid oxide electrolysis cell, the heat exchanger network and the heat supply of a PtL-plant. PtL-efficiency and production costs are used as objectives to generate multiple Pareto fronts for future price scenarios. The results show that the sensitivity to price changes has different impacts on design and operating parameters, which can lead to unattractive solution domains in the Pareto front. Currently, synthetic fuels can be produced at $\numrange{1.83}{2.36}$ $\SI{}{\EUR\per\kilogram}$. In the best case, at $\numrange{1.42}{1.97}$ $\SI{}{\EUR\per\kilogram}$ and $\numrange{3.88}{4.28}$ $\SI{}{\EUR\per\kilogram}$ in the worst case. This paper supports decision-makers in planning PtL-plants to ensure sustainable synthetic fuel availability on a global scale.
\end{abstract}




\end{frontmatter}


\section{Introduction}
\label{sec:introduction}

The leading causes of anthropogenic climate change are \ch{CO2} emissions from the combustion of fossil fuels. The transport sector contributed with \SI{7.95}{\giga\tonne} to approx. \SI{21.60}{\percent} of the annual \ch{CO2} emissions in 2022 \cite{international_energy_agency_co2_2022, international_energy_agency_global_2022}. Facing these figures, there is an increasing need to create a climate-neutral transport sector. Synthetic fuels from renewable energy sources and \ch{CO2} offer a promising solution. They significantly reduce the carbon footprint from the transport sector since the combustion process only releases the \ch{CO2} that was previously taken out of the atmosphere or has been emitted from another source, such as a cement plant \cite{ueckerdt_potential_2021}. Unlike fossil fuels, no previously bonded \ch{CO2} is released into the atmosphere. Promising applications are conceivable both in the short term for existing passenger cars and in the long term for non-electrifiable sectors such as shipping, aviation and heavy machinery.

Today's Power-to-liquid (PtL) plants already produce several thousand tons of methanol and Fischer-Tropsch (FT) fuels annually. A list of currently operating plants and upcoming projects is provided by Pratschner et al. \cite{pratschner_evaluation_2023}. There is also an interactive map with additional sites in \cite{efuel_alliance_ev_efuel_2023}. Another interesting map is the Power-to-X potential atlas by Fraunhofer IEE \cite{fraunhofer_iee_global_2023}, which illustrates the enormous potential of synthetic fuel production sites. So far, little is known about the production costs of the existing plants. The demonstration plant in Haru Oni, Chile, for example, produces fuel at costs of \SI{50}{\EUR\per\liter} \cite{siemens_energy_haru_2023}. However, these costs are not representative for large-scale industrial production. They nevertheless reflect the problem of economically viable production. The production costs must be lowered to be financially relevant to end-users \cite{ngando_co_2022}. The analyses of Ueckerdt et al. \cite{ueckerdt_potential_2021} predict long-term production costs of less than \SI{1}{\EUR\per\liter}. With production costs in the same order of magnitude as conventional fuels, synthetic fuels will be deployed. 

The crucial challenges include producing fuels with competitive costs and in sufficient quantities. The production of synthetic fuels requires a complex interaction of processes, where various parameters, such as the operating point and heat integration, significantly impact costs and efficiency. Applying mathematical programming, the plant design, respectively, the interaction of sub-components can be optimized to meet defined objectives. However, optimizing only a constrained single objective problem may result in insufficient performance of other critical aspects compared to multi-criteria optimization (MOO) \cite{mahrach_comparison_2020}. In this paper, the objectives PtL-efficiency and production costs are used.

The electricity price is another crucial factor that significantly impacts the production costs of synthetic fuels. Although production costs are not directly affected by \ch{CO2} prices, the fuel price for end users is nevertheless affected. For the end users, the point at which the \ch{CO2} prices are applied is not decisive since they are added to the fuel costs anyway. In this paper, we consider the \ch{CO2} price as part of the production cost to ensure comparability with conventionally produced fuels. The prices for renewable electricity and \ch{CO2} certificates can be subject to significant fluctuations due to technological leaps and changing policies. These fluctuations can have unforeseen effects on the economics of synthetic fuel costs. For the design optimization of the PtL-plant, defined assumptions must be made for electricity and \ch{CO2}  prices. Changing prices after the commissioning of the plant can have unwanted effects on production costs. During the planning phase, it is essential to analyze the effects of price changes and consider their consequences when choosing the optimal plant design and operation.

\subsection{Novelty \& Contribution}
In this paper, we conduct a parametric study to understand the impact of changing electricity and \ch{CO2} prices on the production costs of synthetic fuels. We perform an in-depth analysis based on the coupled optimization of the operating characteristics, the heat exchanger network (HEN) and the internal heat supply. PtL-efficiency and fuel production cost are used as objectives of the MOO problem. Starting from a Pareto front at current feedstock prices, the \ch{CO2} and electricity prices are varied. We examine each cost parameter's sensitivity to operating and design characteristics and derive projected Pareto fronts within reasonable price scenarios. We show that fluctuating cost parameters affect design and operation parameters differently. Accordingly, the relevant region of the Pareto front can be narrowed and the plant can be designed to be more resilient. Our analyses are crucial to assess fuel production costs to changing market conditions and policies. This allows us to ensure that synthetic fuels can be made economically viable and available in sufficient quantities to meet global demand.

\subsection{Paper Organization}
The novel \SI{1}{\mega\watt} PtL-plant with its main components and characteristics is presented in Section \ref{sec:systemDescription}. In Section \ref{sec:methods}, the methods for optimization, linearization and the transfer to MILP are presented. In Section \ref{sec:objectives} the modeling and the two antagonistic objectives, PtL-efficiency and production cost, are presented. Further, in Section \ref{sec:CostParameters} the correlation of the cost parameters with the objective functions and the parameter domain are described. In Section \ref{sec:reults}, the parameter study results and the cost parameters' sensitivities are analyzed and discussed.
\section{Materials \& Methods}
\label{sec:MaterialsMethods}

\subsection{System Description}
\label{sec:systemDescription}
\begin{figure*}[htp]
    \centering
    \includegraphics[width=1\textwidth]{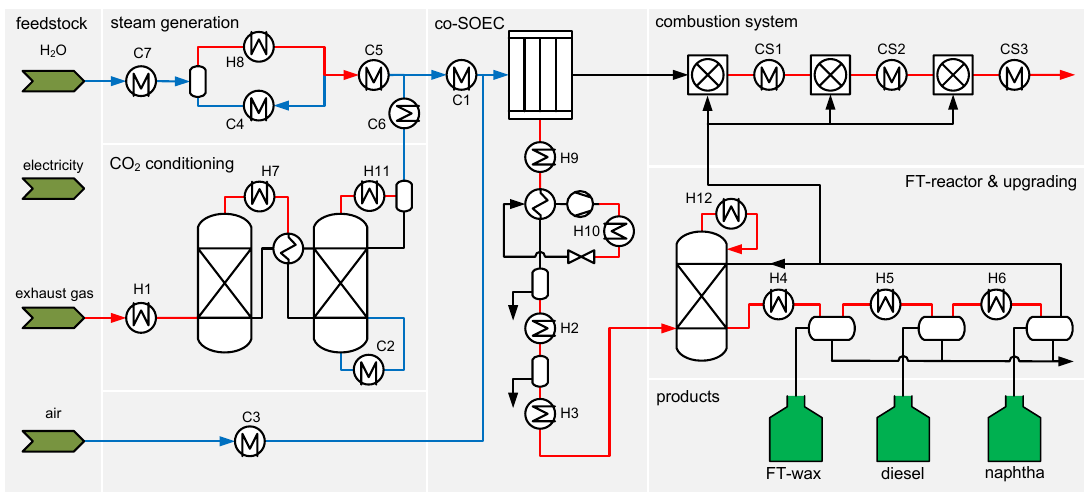}
    \caption{Schematic representation of the \mbox{$\SI{1}{\mega\watt}$} PtL-plant with the five main components and the heat exchangers. Adapted from \cite{huber_unlocking_2023}.} \label{fig:schematicPfd}
\end{figure*}

In the context of the \textit{IFE} (de.: Innovation Flüssige Energie, eng.: Innovation Liquid Energy) research project, a PtL-plant is designed to produce synthetic fuels using water, renewable electricity, exhaust gas from a cement plant and air as feedstock. The plant is designed with a maximum electrolysis capacity of approximately \mbox{\SI{1}{\mega\watt}} using a high-temperature solid oxide co-electrolysis (co-SOEC). Figure \ref{fig:schematicPfd} shows a schematic of the process. Notably, the schematic does not include valves, pumps, or compressors but indicates heat exchangers (HEX) for heat transfer. The cold process streams are shown in blue, and the hot process streams are shown in red.

The PtL process is divided into the following five main components. A more detailed description of the process can be found in \cite{huber_unlocking_2023, huber_ECOS_HENSIFEusecase}.

\begin{description}
    \item[Steam Generation] Pure water at 20~°C is preheated, evaporated, and superheated in the steam generator. 
    \item[\ch{CO2} Conditioning] The PtL-process utilizes exhaust gas from a cement plant with \SI{15.17}{\wtpercent} \ch{CO2}, \SI{81.11}{\wtpercent} \ch{N2}, and \SI{3.72}{\wtpercent} \ch{H2O} at \SI{40}{\celsius}. The \ch{CO2} is conditioned to a purity of \SI{98.73}{\wtpercent} with low residual water content through adsorption and desorption processes.
    \item[co-SOEC] The central component is the co-SOEC, where conditioned \ch{CO2} is superheated, mixed with superheated steam and reformed with preheated air to an \ch{H2}-rich gas and \ch{CO}. The synthesis gas leaving the co-SOEC is cooled and condensed in four stages before entering the FT-reactor.
    \item[FT-reactor and upgrading] A catalytic conversion of the synthesis gas at high temperatures and pressure is carried out in the FT-reactor. The resulting syncrude is separated into FT-wax, diesel and naphtha and upgraded for final use. The product properties downstream of the separation are given in \ref{sec:productprop}. Unreacted synthesis gas is partially recirculated and fed into the combustion system. 
    \item[Combustion System] The combustion system (CS) comprises three serially connected combustion chambers. Unreacted synthesis gas from the separation is used as fuel gas. The CS serves as an internal hot utility to heat the cold process streams.
\end{description}

\subsection{Methods}
\label{sec:methods}
The operating point of the co-SOEC influences the stream parameters inlet, outlet temperature and heat capacity flow of the HEN and vice versa; a holistic optimum can only be achieved by coupling operation and design optimization. The optimization performed in this paper is based on the extension of the classical heat exchanger network synthesis (HENS) by Huber et al. \cite{huber_unlocking_2023, huber_hens_2023}. The HENS of Yee \& Grossman \cite{yee_simultaneous_1990} has been adapted to implement streams with variable temperatures and heat flow capacities. Furthermore, the formulation has been extended to include operating characteristics in the optimization problem. A fundamental assumption for this is that distinct operating points affect only the stream parameters and the objectives. For each operating variable, a piece-wise linear model is created to represent the interaction between the operating point and the stream parameters. For further information regarding the coupled optimization, we refer to \cite{huber_unlocking_2023}.

In the case of the PtL-plant, the costs for \ch{CO2} and electricity significantly influence the production costs. Therefore, the multi-criteria optimization problem is solved in the first instance, considering today's electricity and \ch{CO2} costs. Based on the resulting non-dominated solutions, the cost parameters are varied to study their impact on PtL-efficiency and production costs. The domain in which the parameters vary is specified, considering future development scenarios with rising and falling prices. In a subsequent step, the objectives are recalculated with modified cost parameters for the given design and operating points of the non-dominated solutions. If objectives are recalculated with different parameters while the variables remain the same, this usually does not lead to optimal results. The choice of parameters affects the solution space of the optimization problem, so optimality cannot be guaranteed. In any case, optimality can be achieved by resolving the optimization problem. However, the objective functions cost parameters we vary do not require re-solving the optimization problems. This remarkably efficient approach is possible because the cost parameters do not directly affect the PtL-efficiency and only linearly affect the production costs. As a result, the non-dominated solutions of the Pareto front are only linearly shifted by the production costs at constant efficiency. This procedure enables us to highly efficiently assess the influence of different price scenarios on the system performance and production costs without performing time-consuming optimizations.

\subsubsection{Multi-Criteria Optimization}
\label{sec:MOO}
In this paper, we use the epsilon constraint method to obtain uniformly distributed solutions on the Pareto front. However, the conventional epsilon constraint method is not able to find solutions in overhanging regions of the Pareto front. To overcome this shortcoming, a double-sided epsilon constraint method \cite{huber_hens_2023} is used. Therefore, both sides of the Pareto front are constrained to force the objective into a defined domain.

\subsubsection{Linearization \& Transfer to MILP}
\label{sec:linearization}
Within the adapted HENS, piece-wise linear approximations are used to model all non-linearities. Both piece-wise convex combinations and plane simplices are used to model the multi-dimensional correlations for HEX areas, energy balances and objectives.

To expedite computation, piece-wise linear approximated functions are transformed into mixed integer linear programming (MILP) with a minimal amount of binary variables. One-dimensional, mainly convex curved functions are converted to MILP without binary variables, while other functions necessitate binary variables. Employing the logarithmic coding approach, as proposed by Vielma and Nemhauser \cite{vielma_modeling_2011}, minimizes the number of binary variables. 

Additional information on the methods applied regarding the piece-wise linear approximation and the conversion to MILP can be found in \cite{huber_hens_2023}.

\section{Modeling}
\label{sec:modeling}
The modeling of the system is based on steady-state simulations with Aspen HYSIS. The main operating parameter of the PtL-plant is the cell voltage of the co-SOEC. The system was simulated for seven equidistant cell voltages between \mbox{$U_{\textrm{cell}}^{\textrm{min}} = \SI{1.275}{\volt}$} and \mbox{$U_{\textrm{cell}}^{\textrm{max}} = \SI{1.305}{\volt}$}. The simulation data is used to model the system. Feedstock, power consumption of the co-SOEC, product output, and stream data for the HENS are modeled as a function of cell voltage. The systems’ subcomponents’ size is independent of the cell voltage, resulting in identical system costs. The size of the CS providing the internal heat is also independent of the cell voltage. Only the available amount of FT-offgas is limited depending on the cell voltage.

Detailed information on system modeling can be found in \cite{huber_unlocking_2023}. The parameters used to specify the streams and the heat exchanger network are given in \ref{sec:henparam}.

\subsection{Objectives}
\label{sec:objectives}
Multi-objective optimization enables a holistic performance evaluation of the PtL-plant. In this paper, the antagonistic objectives PtL-efficiency and fuel production costs are optimized.

\subsubsection{PtL-Efficiency}
\label{sec:objPtL}
PtL-efficiency is maximized and described as the ratio of chemically bounded energy in the product \mbox{$\dot H_{\mathrm{prod}}$} to electrical energy input \mbox{$P_{\textrm{el}}$} according to Equation \eqref{eq:objEta}.
\begin{equation}
    \max \eta_{\textrm{PtL}} = \frac{\dot H_{\mathrm{prod}}}{P_{\textrm{el}}} =
    \frac{ \sum_{v}^{} \dot{m}_{\mathrm{prod,}v} \, h_{\mathrm{prod,}v}}{
    P_{\textrm{sys}} + \sum_{j} \varepsilon_{\textrm{hu}} \, q_{\mathrm{hu},j} + \sum_{i} \varepsilon_{\textrm{cu}} \, q_{\mathrm{cu},i}}
\label{eq:objEta}
\end{equation}
The chemically bonded energy downstream of the separation is derived from the sum of the product mass flow rates and the specific enthalpies.

A significant part of the total electrical energy demand $P_{\textrm{el}}$ is the system power $P_{\textrm{sys}}$. With $P_{\textrm{sys}}$, the electric demand of the co-SOEC, circulation pumps, valves and control equipment and losses are covered. The two sums in the denominator represent the energy demand of the electrified utilities. A coefficient of performance of $\varepsilon_{\textrm{cu}} = 0.05$ is assumed for the cold utilities. For the hot utilities, $\varepsilon_{\textrm{hu}}$ is assumed to be 1.05.

\subsubsection{Fuel Production Costs}
\label{sec:objcProd}
The fuel production costs \mbox{$c_{\textrm{prod}}$} are minimized and modeled according to Equation \eqref{eq:objFuel} as the ratio of total annual costs $\mathit{TAC}$ to total synthetic fuel output \mbox{$\sum_{v}^{} t \, \dot{m}_{\mathrm{prod,}v}$}.

\begin{equation}
    \min c_{\textrm{prod}} = \frac{ \mathit{TAC}}{\sum_{v}^{} t \, \dot{m}_{\mathrm{prod,}v}} = \frac{\mathit{CAPEX} + \mathit{OPEX}}{\sum_{v}^{} t \, \dot{m}_{\mathrm{prod,}v}}
\label{eq:objFuel}
\end{equation}
The capital expenditures $\mathit{CAPEX}$ according to Equation \eqref{eq:CAPEX} are composed of investment costs for the system \mbox{$C_{\textrm{sys}}$} and costs for the heat exchanger network according to Yee \& Grossmann \cite{yee_simultaneous_1990}.
\begin{equation}
\begin{split}
    &\mathit{CAPEX} = \\ 
    &\mathit{AF}_{\mathrm{inv}} \left[ \underbrace{C_{\textrm{sys}}}_\text{investment costs}  + \underbrace{\sum_{i}\sum_{j}\sum_{k}{c_\mathrm{v,hex}\left(\frac{q_{ijk}}{U_{ij}\, \mathit{LMTD}_{ijk}} \right)^{\beta}}}_\text{variable HEX stream costs} \right. \\
    &+ \left.\underbrace{\sum_{i}{c_\mathrm{v,hex} \left(\frac{q_{\mathrm{cu},i}}{U_{\mathrm{cu},i}\, \mathit{LMTD}_{\mathrm{cu},i}} \right)^{\beta}}}_\text{variable HEX cold utility costs} \right.\\
    &+ \left. \underbrace{\sum_{j}{c_\mathrm{v,hex} \left(\frac{q_{\mathrm{hu},j}}{U_{\mathrm{hu},j}\, \mathit{LMTD}_{\mathrm{hu},j}} \right)^{\beta}}}_\text{variable HEX hot utility costs} \right. \\
    &+ \left. \underbrace{\sum_{i}\sum_{j}\sum_{k}{c_{\mathrm{f,hex}}\,z_{ijk}} + \sum_{i}{c_{\mathrm{f,hex}} \, z_{\mathrm{cu},i}} + \sum_{j}{c_{\mathrm{f,hex}} \, z_{\mathrm{hu},j}}}_\text{fixed investment costs hex} \right]
\label{eq:CAPEX}
\end{split}
\end{equation}
The operating expenditures $\mathit{OPEX}$ according to Equation \eqref{eq:OPEX} are feedstock and electricity costs depending on the annual full load hours $t$.
\begin{equation}
\begin{split}
    &\mathit{OPEX} = \\
    & \mathit{AF}_{\mathrm{op}} \, t \left[ 
    \underbrace{c_{\textrm{\ch{CO2}}} \, \dot m_{\textrm{\ch{CO2}}} + c_{\textrm{\ch{H2O}}} \, \dot m_{\textrm{\ch{H2O}}} + c_{\textrm{air}} \, \dot m_{\textrm{\ch{air}}}}_\text{feedstock costs}  \right. \\
    &+ \left. \underbrace{ c_{\textrm{el}} \, \left(  P_{\textrm{sys}} +  \sum_{j} \varepsilon_{\textrm{uh}} \, q_{\mathrm{uh},j} +  \sum_{i} \varepsilon_{\textrm{uc}} \, q_{\mathrm{uc},i} \right)}_\text{electricity costs}  \right]
\label{eq:OPEX}
\end{split}
\end{equation}

\subsection{Cost Parameters}
\label{sec:CostParameters}
\begin{figure}[htp]
    \centering
    \includegraphics[width=1\linewidth]{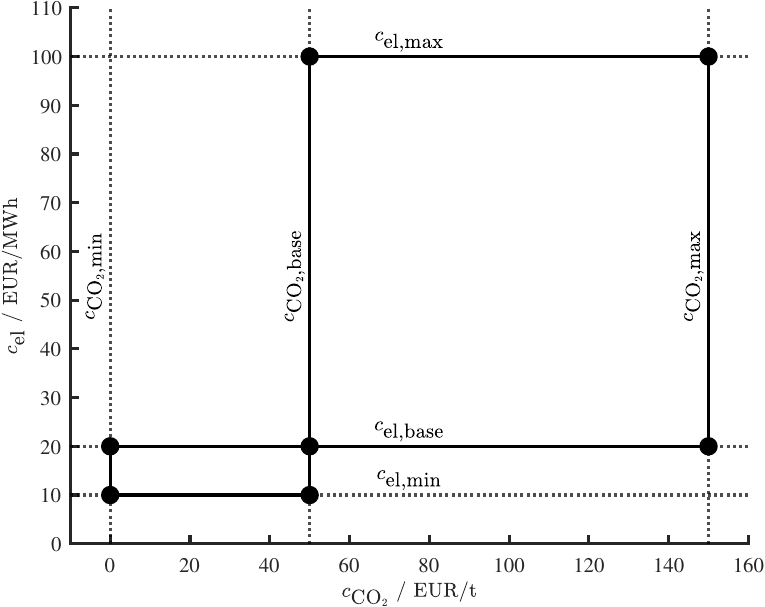}
    \caption{Upper and lower limits of the cost parameters for electricity $c_{\textrm{el}}$ and \ch{CO2} $c_{\textrm{\ch{CO2}}}$.}
    \label{fig:paramRange}
\end{figure}

The coupled optimization of the PtL-plant is performed with initial values for electricity and \ch{CO2} costs based on current market conditions. Figure~\ref{fig:paramRange} shows the chosen domain of cost parameters. The production of synthetic fuels can only be climate neutral if only \ch{CO2}-neutral produced electricity is used. For a location in central Europe and electricity production from wind and solar PV, Janssen et al. \cite{janssen_country-specific_2022} suggest a price of \mbox{$c_{\textrm{el,base}}=\SI{20}{\EUR\per\mega\watt\per\hour}$}. For the lower electricity price limit \mbox{$c_{\textrm{el,min}}$}, we refer to Sens et al. \cite{sens_capital_2022} where in 2050, electricity from PV and wind onshore can be expected to have a levelized cost of electricity (LCOE) of approximately \mbox{$\SI{10}{\EUR\per\mega\watt\per\hour}$}. We consider this as a best-case scenario for future electricity price developments. The upper limit of \mbox{$c_{\textrm{el,max}}=\SI{100}{\EUR\per\mega\watt\per\hour}$} is derived from the LCOE for wind offshore in 2050 \cite{sens_capital_2022}.

The price for emitted \ch{CO2} can be levied either with a cap-and-trade system like the Emission Trading System (ETS) or as a carbon tax \cite{kruse-andersen_optimal_2022}. Both policies aim to have the government financially penalize greenhouse gas emissions and thus force the polluters to reduce. What \ch{CO2} pricing will look like for synthetic fuels has not yet been defined. In any case, it can be assumed that the end users will bear the costs. In this paper, we use a \ch{CO2} price of \mbox{$c_{\textrm{\ch{CO2},base}}=\SI{50}{\EUR\per\tonne}$} as a base value. This corresponds roughly to the average \ch{CO2} tax in Europe according to the Carbon Pricing Dashboard of The World Bank \cite{the_world_bank_carbon_2023}. A minimum price of \mbox{$c_{\textrm{\ch{CO2},min}}=\SI{0}{\EUR\per\tonne}$} is assumed as an economic best-case scenario. This scenario represents sites with no \ch{CO2} pricing or where special regulations for synthetic fuels have been enacted. As a maximum value, a \ch{CO2} price of \mbox{$c_{\textrm{\ch{CO2},min}}=\SI{150}{\EUR\per\tonne}$} is assumed. This represents a roadmap for realistic price development towards 2050 \cite{bloombergnef_carbon_2023, noauthor_essential_2022, noauthor_ghg_2022}.

\subsection{Implementation}
\label{sec:implementation}

The optimization problem was formulated with Yalmip R20210331 \cite{lofberg_toolbox_2004} and Matlab R2022b \cite{the_mathworks_inc_matlab_2022}. Gurobi 10.0.0 was used as MILP solver \cite{gurobi_optimization_llc_gurobi_2023}. A MIP gap of less than \SI{1}{\percent} was set as a termination criterion for the optimization.  All computations were performed on a 64-core server (AMD EPYC 7702P) with 265 GB of RAM.

All non-linearities are piecewise linearly approximated, analogous to the approach from Huber et al. \cite{huber_unlocking_2023}. A root-mean-square error (RMSE) of less than \SI{1}{\percent} is set as a termination criterion for the refinement of the approximations.

All cost parameters are given in \ref{sec:costParam}.
\section{Results}
\label{sec:reults}

\begin{figure}[H]
    \centering
    \includegraphics[width=1\linewidth]{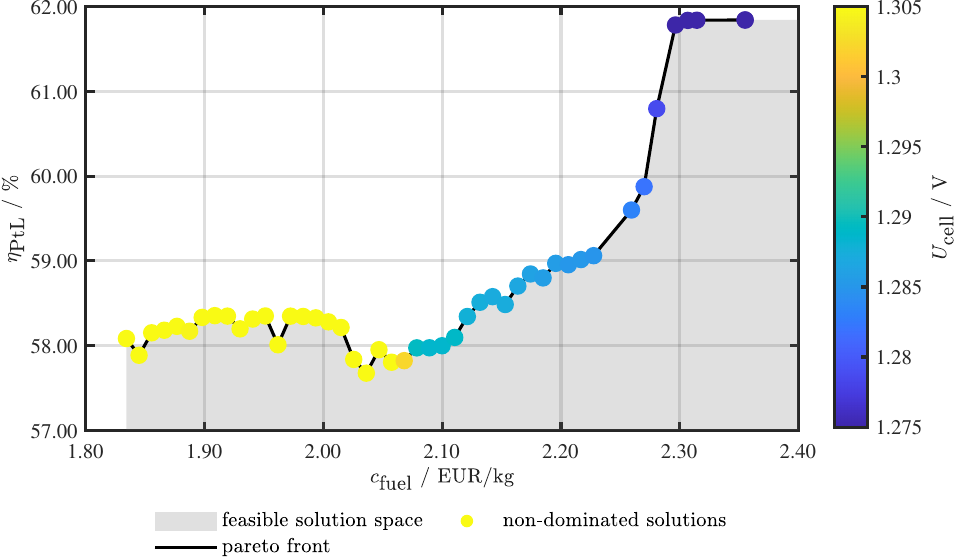}
    \caption{Non-dominated solutions and Pareto front of the coupled optimization. The color coding represents the cell voltage as the main operating parameter.}
    \label{fig:paretoInitial}
\end{figure}

The non-dominated solutions of the antagonistic objectives are illustrated as a Pareto front in Figure~\ref{fig:paretoInitial}. The color coding represents the cell voltage as the critical operating parameter. The slight scattering of the Pareto front at low production costs results from the solver time out and the optimality gap of the solution. There were no solutions found for the gaps around high production costs within the defined solver time. The stream plots at the corner points of the Pareto front and the characteristic values such as number the of heat exchangers, power of the utilities and the costs can be obtained from \cite{huber_unlocking_2023}.

An essential aspect of the Pareto front is that the production costs increase with decreasing cell voltage. Further, the shape of the Pareto front allows us to derive the impact of design and operating parameters. At the highest cell voltage of $U_{\textrm{cell}} = \SI{1.305}{\volt}$, the PtL-efficiency remains almost unchanged between $\SI{57.67}{\percent}$ and $\SI{58.35}{\percent}$. The production costs, however, vary from \mbox{$\SI{1.83}{\EUR\per\kilogram}$} to \mbox{$\SI{2.06}{\EUR\per\kilogram}$}. This indicates that design parameters significantly influence production costs but only bearly the efficiency in this region. Conversely, at the upper production cost range, the influence of the cell voltage dominates. The production costs remain almost constant while the PtL-efficiency increases significantly. 

Price changes are analyzed based on the non-dominated solutions optimized with initial cost parameters. Production costs are recalculated for each solution at constant PtL-efficiencies. This is feasible because the cost parameters for \ch{CO2} and electricity occur as linear terms only in the objective function of the production costs and not in that of the PtL-efficiency; see Equation \eqref{eq:objFuel} and \eqref{eq:OPEX}, respectively. Thus, only the production costs change depending on the cost parameters, while the PtL-efficiency remains unchanged. We quantify the impact of the cost parameters with the help of the parameter sensitivity defined by Equation \eqref{eq:sencel} and \eqref{eq:sencCO2}. 
\begin{equation}
\label{eq:sencel}
    \nabla c_{\textrm{el}} = \mathit{AF}_{\mathrm{op}} \frac{ \left[ P_{\textrm{sys}} +  \left( \sum_{j} \varepsilon_{\textrm{uh}} \, q_{\mathrm{uh},j} +  \sum_{i} \varepsilon_{\textrm{uc}} \, q_{\mathrm{uc},i} \right) \right]}{\sum_{v}^{} \dot{m}_{\mathrm{prod,}v}}
\end{equation}
\begin{equation}
\label{eq:sencCO2}
    \nabla c_{\textrm{\ch{CO2}}} = \mathit{AF}_{\mathrm{op}}\frac{\dot m_{\textrm{\ch{CO2}}}}{\sum_{v}^{} \dot{m}_{\mathrm{prod,}v}}
\end{equation}

\begin{figure}[htp]
    \centering
    \includegraphics[width=1\linewidth]{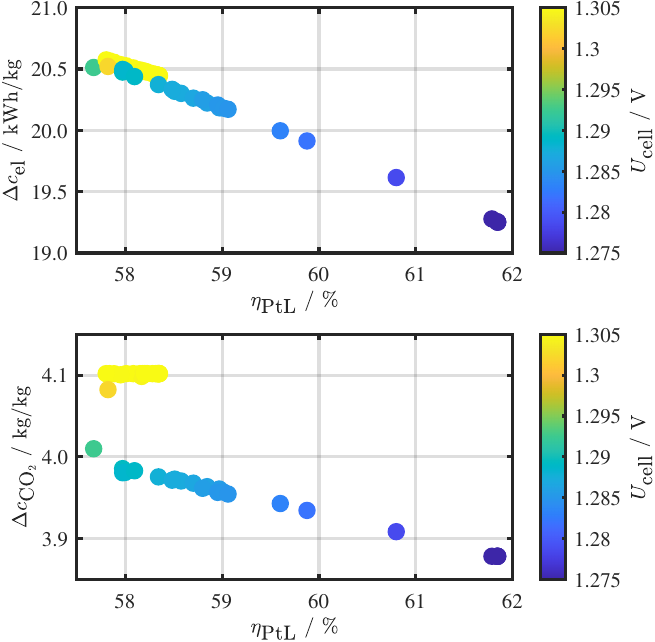}
    \caption{Calculated sensitivities for electricity and \ch{CO2} price changes as a function of PtL-efficiency.}
    \label{fig:sensitivity}
\end{figure}

Figure~\ref{fig:sensitivity} on the top shows the sensitivities to electricity price changes. With an annualization factor of $\mathit{AF}_{\mathrm{op}} = \SI{1}{}$, $\nabla c_{\textrm{el}}$ can be interpreted as energy consumption to produce $\SI{1}{\kilogram}$ of synthetic fuel, respectively FT-wax, diesel and naptha. Fasihi et al. obtained a slightly higher value of $\SI{21.98}{\kilo\watt\hour\per\kilogram}$ for a similar PtL-plant with an ambient air scrubber, reverse water gas shift reaction (RWGS) and alkaline electrolysis cell \cite{fasihi_techno-economic_2016}. In our case, PtL-efficiency strongly correlates with the sensitivity. Higher efficiencies are, therefore, less susceptible to price changes.

The sensitivities to \ch{CO2} price changes are shown in Figure~\ref{fig:sensitivity} on the bottom. They can be interpreted as the utilized amount of \ch{CO2} per kg of synthetic fuel produced. The process presented by Fasani et al. \cite{fasihi_techno-economic_2016} requires $\SI{3,20}{\kilogram\per\kilogram}$. In this case, the sensitivity correlates less strongly with the efficiency but instead much more with the cell voltage. Therefore, plant designs with low cell voltages are less susceptible to price changes.

\begin{figure*}[htp]
    \centering
    \includegraphics[width=1\textwidth]{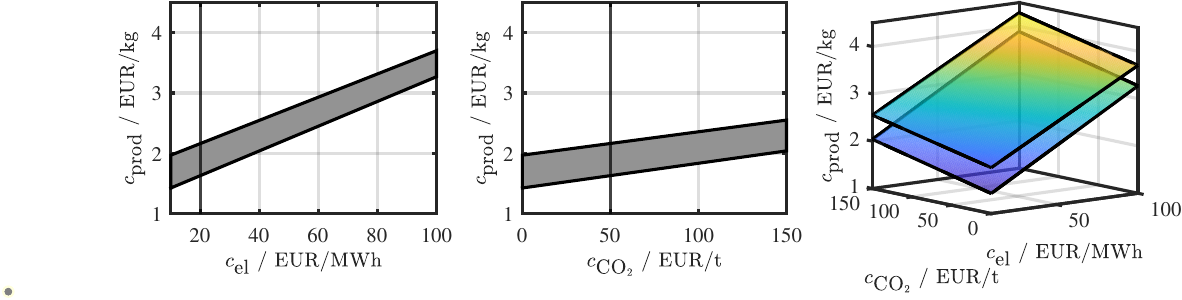}
    \caption{Upper and lower production costs as a function of electricity and \ch{CO2} price. Left: Isolated impact of electricity price. Middle: Isolated impact of the \ch{CO2} price. Right: Combined impact.}
    \label{fig:surface}
\end{figure*}

Based on the sensitivities, we calculated the production costs for each non-dominated solution as function of the electricity and \ch{CO2} prices. Figure~\ref{fig:surface} shows the highest and lowest production costs. The minimum and maximum production costs are shown as a function of electricity prices on the left and \ch{CO2} prices on the right. The horizontal lines show the current cost parameters \mbox{$c_{\textrm{el,base}}$} and \mbox{$c_{\textrm{\ch{CO2},base}}$} from Figure~\ref{fig:paramRange}. The minimum and maximum production costs are shown for both cost parameters on the right side of Figure~\ref{fig:surface}. All non-dominated solutions are between these two planes. In this context, the sensitivities can be interpreted as slopes of the limiting lines, respectively, planes. Accordingly, the higher the sensitivity, the more sensitive the solution reacts to price changes. It can be concluded from the plane slops that the influence of the electricity price on the production costs is about five times larger than that of the \ch{CO2} price.

\begin{figure}[htp]
    \centering
    \includegraphics[width=1\linewidth]{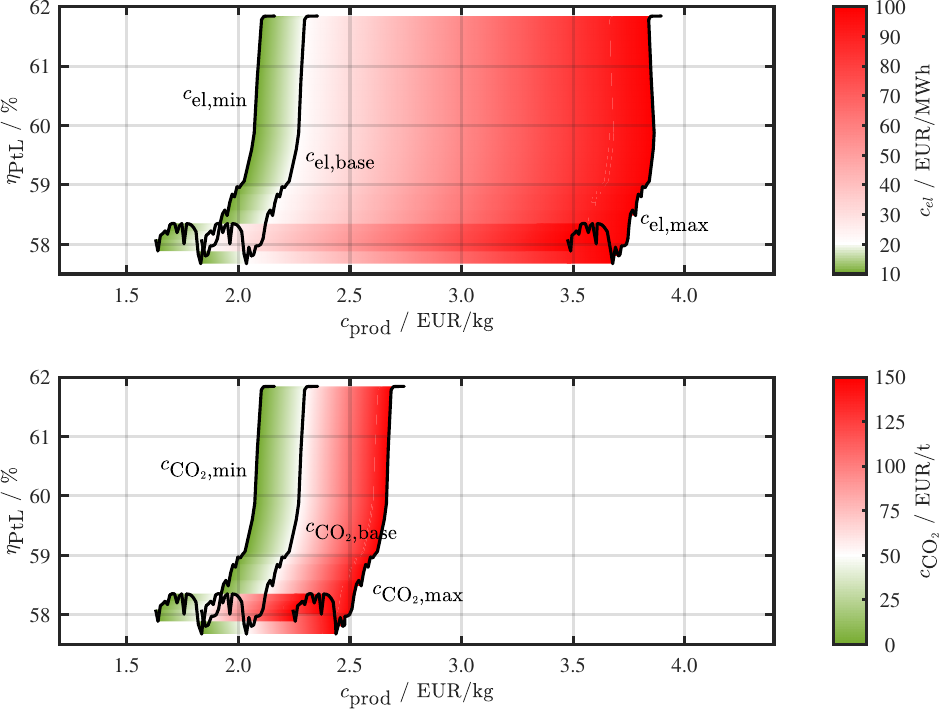}
    \caption{Pareto fronts for changing electricity prices (top) and \ch{CO2} prices (bottom). The color gradient represents the shift of the Pareto front due to price changes.}
    \label{fig:pareto_minmax}
\end{figure}

The price changes affect the shape of the Pareto front depending on the efficiency. On the top of Figure~\ref{fig:pareto_minmax}, three Pareto fronts for minimum, initial and maximum electricity prices are shown. Electricity prices have a substantial effect on production costs. When electricity costs are high, it is noticeable that the Pareto front forms a pocket in the efficiency range from \mbox{$\SI{59.0}{\percent}$} to \mbox{$\SI{61.8}{\percent}$}. If higher electricity costs are expected, selecting a plant design with these parameters should be avoided since lower production costs can be expected at higher and lower efficiencies. The Pareto fronts for different \ch{CO2} prices in Figure~\ref{fig:pareto_minmax} on the bottom do not create pockets. It is noticeable that the influence on the production costs is significantly lower.

\begin{figure}[htp]
    \centering
    \includegraphics[width=1\linewidth]{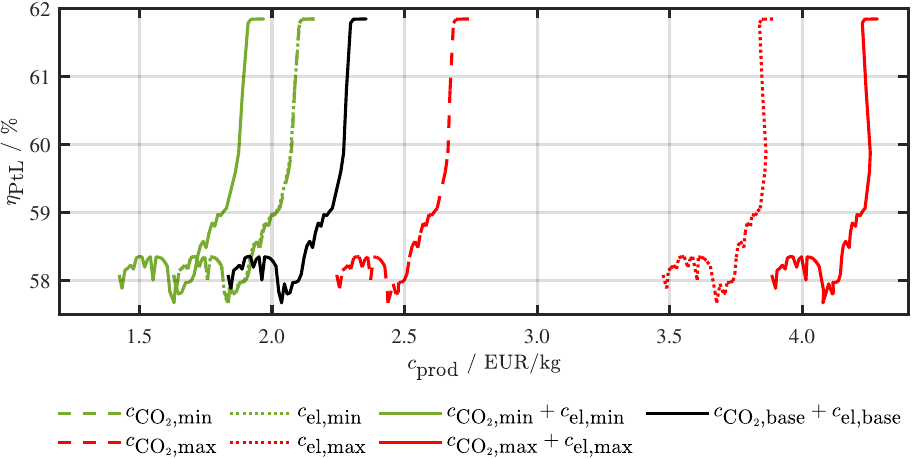}
    \caption{Pareto fronts for the extreme values of the cost parameters: Initial costs (black), cost reduction (green), cost increase (red).}
    \label{fig:pareto_superPos}
\end{figure}

Figure~\ref{fig:pareto_superPos} shows the collective impact of price changes. For each extreme value, a Pareto front is shown. It should be noted that the lines for \mbox{$c_{\ch{CO2},\textrm{min}}$} and \mbox{$c_{\ch{el},\textrm{min}}$} overlap almost completely. If the electricity price is reduced to \mbox{$\SI{10}{\EUR\per\mega\watt}$} and the \ch{CO2} price is dropped to $0$, synthetic fuels can be produced in the \mbox{$\SI{1.42}{\EUR\per\kilogram}$} to \mbox{$\SI{1.97}{\EUR\per\kilogram}$} range. Compared to the current market price for gasoline of about \mbox{$\SI{0.5}{\EUR\per\kilogram}$}, the large-scale introduction of synthetic fuels will be an economic challenge \cite{us_energy_information_administration_us_2023}. At the highest electricity and \ch{CO2} prices of \mbox{$\SI{100}{\EUR\per\mega\watt}$} and \mbox{$\SI{150}{\EUR\per\tonne}$}, the production costs range from \mbox{$\SI{3.88}{\EUR\per\kilogram}$} to \mbox{$\SI{4.28}{\EUR\per\kilogram}$}. Whether synthetic fuels will prevail at this price level is questionable.

\section{Conclusion}
\label{sec:conclusion}
In this paper, we investigated the impact of changing electricity and \ch{CO2} prices on the production costs of synthetic fuels. A novel \mbox{$\SI{1}{\mega\watt}$} PtL-plant was considered as a use case. The main components were modeled and the cell voltage of the co-SOEC was optimized as a crucial operating parameter coupled with the heat exchanger network and the internal heat supply as design parameters. A formulation for PtL-efficiency and production cost was presented as objective functions. Based on the optimization results, the sensitivities for electricity and \ch{CO2} price changes were calculated and their correlation with the design and operational parameters was discussed. A comprehensive cost analysis was performed with defined scenarios for decreasing and increasing cost parameters.

Our results show that the non-dominated solutions of the Pareto front in the low-efficiency area depend more on the design parameters than on the cell voltage of the co-SOEC as the operating parameter. On the other hand, a distinct dependence on the cell voltage can be seen in the high-efficiency area. Current electrictiy and \ch{CO2} prices results in production costs in the range of $\numrange{1.83}{2.36}$ $\SI{}{\EUR\per\kilogram}$ and a PtL-efficiency of $\numrange{57.67}{61.84}$ $\SI{}{\percent}$. Further, we were able to show that the influence of the electricity price is about five times larger than that of the \ch{CO2} price. The sensitivity of the electricity price correlates strongly inversely with the efficiency. This implies that higher efficiencies are less affected by price increases. Substantial price increases lead to pockets in the Pareto front, which may preclude specific design and operating areas. A further key outcome is that our results allow us to estimate the expected production costs depending on price changes. In the economic best-case scenario, production costs can be reduced to the range of $\numrange{1.42}{1.97}$ $\SI{}{\EUR\per\kilogram}$. In the worst-case scenario, production costs increase up to $\numrange{3.88}{4.28}$ $\SI{}{\EUR\per\kilogram}$.

With this paper, we provide an overview of current production costs and the influence of design and operating parameters on future cost developments. A central and, above all, necessary basis for decision-making has been developed with this paper. It enables the decision-makers to evaluate the design and its consequences for sustainable synthetic fuel production. With this paper, we promote the efficient and cost-effective production of synthetic fuels and contribute significantly to a climate-neutral and sustainable future.

\bibliographystyle{elsarticle-num-names} 
\bibliography{references}

\section*{Statements and Declarations}
\subsection*{Funding}
This research was funded by the Austrian Research Promotion Agency (FFG) under grant number 884340 and TU Wien Bibliothek through its Open Access Funding Programme.

\subsection*{Competing Interests}
The authors have no relevant financial or non-financial interests to disclose.

\subsection*{Authors Contributions}
The method presented in this paper was developed by David Huber. Testing and evaluation was done by David Huber. The conceptualization of the paper was done by all authors. The first draft was written by David Huber. All authors contributed to the revision of the initial draft. Funding and supervision was done by René Hofmann. All authors read and approved the final manuscript.

\appendix
\section{Product Properties}
\label{sec:productprop}
Table \ref{tab:productProp} lists the chemical and physical product parameters. These are invariant for the system regardless of the cell voltage.

\begin{table}[H]
    \centering
    \caption{Chemical and physical product properties at $\SI{40}{\celsius}$ and $\SI{101324.97}{\pascal}$ downstream the upgrading. Adapted from \cite{huber_unlocking_2023}.}
    \label{tab:productProp}
    \resizebox{\linewidth}{!}{%
    \begin{tabular}{ccrrr}
        \toprule
        $v$ & product & \multicolumn{1}{c}{$h_{\mathrm{prod}}$ / $\SI{}{\mega\joule\per\kilogram}$} & \multicolumn{1}{c}{$\rho_{\mathrm{prod}}$ / $\SI{}{\kilogram\per\cubic\meter}$} & \multicolumn{1}{c}{$\mu_{\mathrm{prod}}$ / $\SI{}{\milli\pascal\per\second}$} \\
        \midrule
        $1$ & FT-wax & $43.887$ & $797.73$ & $6.7477$ \\
        $2$ & diesel & $44.345$ & $748.81$ & $1.5983$ \\
        $3$ & naphtha & $44.676$ & $516.17$ & $0.5893$ \\
        \bottomrule
    \end{tabular}}
\end{table}

\section{Stream Data}
\label{sec:henparam}

Table~\ref{tab:streamData} provides the stream parameters or their bounds. Values without brackets are independent of the cell voltage constant. Values in square brackets are bounds of stream parameters dependent on cell voltage. The heat transfer coefficients are \mbox{$U = \SI{0.5}{\kilo\watt\per\meter\squared\per\kelvin}$} for all streams.

The HEN problem has been modeled with \mbox{$N_{\textrm{st}}=3$} stages for the heat exchange. The minimum temperature difference of \mbox{$\Delta T_{\textrm{min}}=\SI{1}{\kelvin}$} must not be exceeded.

\begin{table}[H]
    \centering
    \caption{Stream data with limits for inlet, outlet temperature and flow capacity. Adapted from \cite{huber_unlocking_2023}.}
    \label{tab:streamData}
    \resizebox{\linewidth}{!}{%
    \begin{tabular}{crrr}
        \toprule
        Stream & $T^{\mathrm{in}}$ / $\SI{}{\celsius}$ & $T^{\mathrm{out}}$ / $\SI{}{\celsius}$ & $F$ / $\SI{}{\kilo\watt\per\kelvin}$ \\
        \midrule 
        H1 & $40.0$ & $35.0$ & $[1.71, 2.16]$\\
        H2 & $[127.9, 131.1]$ & $[34.0, 35.0]$ & $[0.09, 0.12]$\\
        H3 & $[169.8, 174.1]$ & $[34.0, 35.0]$ & $[0.09, 0.12]$\\
        H4 & $210.0$ & $190.0$ & $[0.27, 0.28]$\\
        H5 & $190.0$ & $120.0$ & $[0.56, 0.58]$\\
        H6 & $120.0$ & $30.0$ & $[0.48, 0.50]$\\
        H7 & $[45.4, 57.0]$ & $31.0$ & $[2.35, 2.95]$\\
        H8 & $138.9$ & $137.9$ & $[59.60, 94.40]$\\
        H9 & $[805.2, 825.5]$ & $[34.0, 35.0]$ & $[0.10, 0.13]$\\
        H10 & $[49.5,50.7]$ & $[34.0, 35.0]$ & $[0.65, 0.88]$\\
        H11 & $101.8$ & $30.0$ & $[0.51 ,0.64]$\\
        H12 & $190.0$ & $188.0$ & $[76.88, 80.45]$\\
        C1 & $[318.0, 319.2]$ & $[825.0, 870.5]$ & $[0.14, 0.18]$\\
        C2 & $116.9$ & $124.2$ & $[20.02, 25.12]$\\
        C3 & $[57.3, 58.8]$ & $825.0$ & $[0.25, 0.33]$\\
        C4 & $137.9$ & $139.9$ & $[105.77, 142.64]$\\
        C5 & $138.9$ & $[426.6, 449.4]$ & $[0.10, 0.11]$\\
        C6 & $35.0$ & $[115.9, 145.4]$ & $[0.05, 0.06]$\\
        C7 & $20.3$ & $[189.5, 199.6]$ & $[0.15, 0.21]$\\
        \midrule
        CS1 & $900.0$ & $[100, 890]$ & $[59.60, 94.40]$\\
        CS2 & $900.0$ & $[100, 890]$ & $[0.10, 0.13]$\\
        CS3 & $900.0$ & $[100, 890]$ & $[0.65, 0.88]$ \\
        \bottomrule
    \end{tabular}}
\end{table}

\section{Cost Parameters}
\label{sec:costParam}
Table \ref{tab:costParam} lists the cost parameters and their sources. The investment costs are depreciated linearly over a period of $\SI{20}{}$ years. Thus, the annual depreciation factor is \mbox{$\mathit{AF}_{\textrm{inv}} = \nicefrac{1}{20 \textrm{y}}$}. The depreciation factor for the operating costs is assumed to be \mbox{$\mathit{AF}_{\textrm{op}} = \SI{1}{\per\year} $}. Analogous to \cite{adelung_global_2022, herz_economic_2021, dieterich_power--liquid_2020}, \mbox{$t = \SI{8000}{\hour\per\year}$} full load operating hours per year are assumed.

\begin{table}[htp]
    \centering
    \caption{Cost parameters and their source. Adapted from \cite{huber_unlocking_2023}.}
    \label{tab:costParam}
    \resizebox{\linewidth}{!}{%
    \begin{tabular}{cll}
        \toprule
        cost share & \multicolumn{1}{c}{value} & \multicolumn{1}{c}{comment / source} \\
        \midrule
        $\beta$ & $\SI{0.8}{}$ & \\
        $c_{\textrm{f,hex}}$ & $\SI{1013.6}{\EUR\per\year}$ & AISI 316, interpolated from \cite{noauthor_dace_2021}\\
        $c_{\textrm{v,hex}}$ & $\SI{61.8}{\EUR\per\betaCo\per\year}$ & \\
        \midrule
        $C_{\textrm{sys}}$ & $\SI{10000000}{\EUR}$ & project internal estimation \& \cite{herz_economic_2021} \\
        \midrule
        $c_{\textrm{\ch{H2O}}}$ & $\SI{3.54}{\EUR\per\tonne}$ & mean for Europe \cite{eureau_europes_2021} \\
        $c_{\textrm{\ch{CO2},base}}$ & $\SI{50}{\EUR\per\tonne}$ & average \ch{CO2} tax in Europe \cite{the_world_bank_carbon_2023} \\
        $c_{\textrm{\ch{air}}}$ & $\SI{0}{\EUR\per\tonne}$ & ambient air is free of charge \\
        $c_{\textrm{\ch{el},base}}$ & $\SI{20}{\EUR\per\mega\watt\per\hour}$ & at a projected plant location in Europe \cite{janssen_country-specific_2022} \\
        \bottomrule
    \end{tabular}}
\end{table}

\makenomenclature

\renewcommand\nomgroup[1]{%
  \item[\bfseries
  \ifstrequal{#1}{A}{Acronyms}{%
  \ifstrequal{#1}{V}{Variables}{%
  \ifstrequal{#1}{L}{Subscripts}{%
  \ifstrequal{#1}{H}{Superscripts}{%
  \ifstrequal{#1}{S}{Sets}}}}}%
]}

\setlength{\nomlabelwidth}{1.5cm}

\newcommand{\nomunit}[1]{%
\renewcommand{\nomentryend}{\hspace*{\fill}#1}}

\nomenclature[A]{MILP}{mixed-integer linear programming}
\nomenclature[A]{MIP}{mixed-integer programming}
\nomenclature[A]{PtL}{power-to-liquid}
\nomenclature[A]{SOEC}{solid oxide electrolysis cell}
\nomenclature[A]{FT}{fischer-tropsch}
\nomenclature[A]{MOO}{multi objective optimization}
\nomenclature[A]{HEN}{heat exchanger network}
\nomenclature[A]{HENS}{heat exchanger network synthesis}
\nomenclature[A]{LCOP}{levelized cost of product}
\nomenclature[A]{IFE}{Innovation Flüssige Energie, eng.: Innovation Liquid Energy}
\nomenclature[A]{CS}{combustion system}
\nomenclature[A]{RMSE}{root-mean-square error}
\nomenclature[A]{CAPEX}{annual capital expenses}
\nomenclature[A]{OPEX}{operational expenditures}
\nomenclature[A]{LMTD}{logarithmic mean temperature difference}
\nomenclature[A]{HEX}{heat exchanger}
\nomenclature[A]{PV}{photovoltaics}
\nomenclature[A]{ETS}{emission trading system}
\nomenclature[A]{RWGS}{reverse water gas shift}

\nomenclature[V]{$F$}{{flow capacity}\nomunit{\SI{}{\kilo\watt\per\kelvin}}}
\nomenclature[V]{$U_{\textrm{cell}}$}{{cell voltage}\nomunit{\SI{}{\volt}}}
\nomenclature[V]{$\dot{m}$}{{mass flow}\nomunit{\SI{}{\kg\per\hour}}}
\nomenclature[V]{$T$}{{temperature}\nomunit{\SI{}{\celsius}}}
\nomenclature[V]{$h_{\textrm{prod}}$}{{specific enthalpy of the product}\nomunit{\SI{}{\mega\joule\per\kilogram}}}
\nomenclature[V]{$\rho$}{{density}\nomunit{\SI{}{\kilogram\cubic\meter}}}
\nomenclature[V]{$\mu$}{{dynamic viscosity}\nomunit{\SI{}{\milli\pascal\per\second}}}
\nomenclature[V]{$N_{\textrm{st}}$}{number of stages}
\nomenclature[V]{$P_{\textrm{sys}}$}{{electrical energy demand w/o utilities}\nomunit{\SI{}{\kW}}}
\nomenclature[V]{$U$}{{overall heat transfer coefficient}\nomunit{\SI{}{\kilo\watt\per\meter\squared\per\kelvin}}}
\nomenclature[V]{$\eta_{\textrm{PtL}}$}{{PtL-efficiency}\nomunit{\SI{}{\percent}}}
\nomenclature[V]{$P_{\textrm{el}}$}{{total electrical energy demand}\nomunit{\SI{}{\kilo\watt}}}
\nomenclature[V]{$\dot{H}$}{{chemically bounded energy in FT-products}\nomunit{\SI{}{\kilo\watt}}}
\nomenclature[V]{$\varepsilon$}{{coefficient of performance}\nomunit{\SI{}{}}}
\nomenclature[V]{$q$}{{heat flow}\nomunit{\SI{}{\kilo\watt}}}
\nomenclature[V]{$t$}{{annual full load hours}\nomunit{\SI{}{\hour\per\year}}}
\nomenclature[V]{$C_{\textrm{sys}}$}{{investment costs}\nomunit{\SI{}{\EUR}}}
\nomenclature[V]{$a$}{{depreciation period}\nomunit{\SI{}{\year}}}
\nomenclature[V]{$AF_{\textrm{inv}}$}{{investment annualization factor}\nomunit{\SI{}{\per\year}}}
\nomenclature[V]{$AF_{\textrm{op}}$}{{operational annualization factor}\nomunit{\SI{}{}}}
\nomenclature[V]{$c_{\textrm{prod}}$}{{product costs}\nomunit{\SI{}{\EUR\per\kilogram}}}
\nomenclature[V]{$c_{\text{f}}$}{{feedstock costs}\nomunit{\SI{}{\EUR\per\tonne}}}
\nomenclature[V]{$c_{\textrm{el}}$}{{electricity costs}\nomunit{\SI{}{\EUR\per\mega\watt\per\hour}}}
\nomenclature[V]{$c_{\textrm{f,hex}}$}{{step-fixed HEX costs}\nomunit{\SI{}{\EUR\per\year}}}
\nomenclature[V]{$c_{\textrm{v,hex}}$}{{variable HEX costs}\nomunit{\SI{}{\EUR\per\betaCo\per\year}}}
\nomenclature[V]{$\beta$}{{cost exponent}\nomunit{\SI{}{}}}
\nomenclature[V]{$z$}{{binary variable for existance of HEX}\nomunit{\SI{}{}}}
\nomenclature[V]{$\Delta T_{\textrm{min}}$}{{minimum temperature difference}\nomunit{\SI{}{\kelvin}}}
\nomenclature[V]{$\mathit{TAC}$}{{total annual costs}\nomunit{\SI{}{\EUR\per\year}}}
\nomenclature[V]{$\mathit{CAPEX}$}{{capital expenditures}\nomunit{\SI{}{\EUR\per\year}}}
\nomenclature[V]{$\mathit{CAPEX}$}{{operating expenditures}\nomunit{\SI{}{\EUR\per\year}}}
\nomenclature[V]{$\mathit{LMTD}$}{{logarithmic mean temperature difference}\nomunit{\SI{}{\kelvin}}}

\nomenclature[L]{$\textrm{prod}$}{production}
\nomenclature[L]{$\textrm{cu}$}{cold utility}
\nomenclature[L]{$\textrm{hu}$}{hot utility}
\nomenclature[L]{$\textrm{hex}$}{heat exchanger}
\nomenclature[L]{$\textrm{el}$}{electric}
\nomenclature[L]{$i$}{hot stream}
\nomenclature[L]{$j$}{cold stream}
\nomenclature[L]{$k$}{stage}

\nomenclature[H]{in}{inlet}
\nomenclature[H]{out}{outlet}
\nomenclature[H]{min}{minimum}
\nomenclature[H]{max}{maximum}

\printnomenclature

\end{document}